\documentclass[letterpaper, 10 pt, conference]{ieeeconf}  % Comment this line out if you need a4paper

\IEEEoverridecommandlockouts                              % This command is only needed if 
\usepackage{url}
\usepackage[utf8x]{inputenc}
\usepackage{standalone}
\usepackage{tikz}
\usepackage{pgfplots}
\tikzstyle{block} = [draw, rectangle, 
    minimum height=3em, minimum width=5em]
\tikzstyle{input} = [coordinate]
\tikzstyle{output} = [coordinate]
\tikzstyle{pinstyle} = [pin edge={to-,thin,black}]
\usetikzlibrary{shapes.geometric,shapes.arrows,decorations.pathmorphing,calc}
\usetikzlibrary{matrix,chains,scopes,positioning,arrows,fit}

\usepackage{subcaption}
\usepackage{siunitx}
\usepackage{epsfig}
\usepackage{times}
\usepackage{amssymb}
\usepackage{dsfont}
\usepackage{multirow}
\usepackage{algorithm}
\usepackage[noend]{algcompatible}

\usepackage{amsmath}
\usepackage{steinmetz}
\usepackage{cite}
\usepackage{graphicx}    % include this line if your document contains figures
\usepackage{amsfonts}

\usepackage{amsthm}
\newtheorem{assumption}{Assumption}
\newtheorem{remark}{Remark}
\pgfplotsset{compat=newest, ticks=none}

\title{\LARGE \bf
Tutoring Reinforcement Learning via Feedback Control}
\author{Francesco De Lellis$^{1}$, Giovanni Russo$^{2,*}$, Mario di Bernardo$^{1,*}$% <-this % stops a space
 \thanks{$^{1}$Department of Electrical Engineering and ICT, University of Naples Federico II, Italy}
 \thanks{$^{2}$ Department of Information and Electrical Engineering and Applied Mathematics, University of Salerno, Italy}
 \thanks{$^{*}$Corresponding authors. {\tt\small mario.dibernardo@unina.it, giovarusso@unisa.it}}}

\begin{document}

\maketitle
\thispagestyle{empty}
\pagestyle{empty}

%%%%%%%%%%%%%%%%%%%%%%%%%%%%%%%%%%%%%%%%%%%%%%%%%%%%%%%%%%%%%%%%%%%%%%%%%%%%%%%%
\begin{abstract}
We introduce a control-tutored reinforcement learning (CTRL) algorithm. The idea is to enhance tabular learning algorithms by means of a control strategy with limited knowledge of the system model. By tutoring the learning process, the learning rate can be substantially reduced. We use the classical problem of stabilizing an inverted pendulum as a benchmark to numerically illustrate the advantages and disadvantages of the approach.
\end{abstract}

%%%%%%%%%%%%%%%%%%%%%%%%%%%%%%%%%%%%%%%%%%%%%%%%%%%%%%%%%%%%%%%%%%%%%%%%%%%%%%%%
\section{Introduction}
Reinforcement learning (RL) \cite{Sutton1998,bertsekas1996neuro} is increasingly used to learn control policies from data \cite{kober2013reinforcement,garcia2015comprehensive,cheng2019end} in a number of different applications. Despite the several advantages of this approach to control, one of its key drawbacks lies is the requirement of performing a typically large number of trials to explore the state-action space and hence learn a sub-optimal control policy for the plant of interest. 
In particular, the control policy is found by exploring  the Markov Decision Process encapsulating the control problem, thus accepting possible failures while learning. Unfortunately, in control applications, long training phases are often unacceptable and failures while learning might lead to unsafe situations. 
Moreover, these applications are often characterized by a continuous state-space, and using RL as is requires a dense discretization of the system state space or a function approximation that is also subject to the learning process. 

To overcome these limitations, many different flavours of RL have been developed in the increasingly vast literature on the problem. First and foremost, model-based reinforcement learning where the introduction of some mathematical model of the plant is used to guarantee some degree of stability of the learning process and partially solve some of the issues mentioned above, e.g. \cite{berkenkamp2017safe, gu2016continuous, deisenroth2011pilco, rosolia2017learning, rathi2020driving}.  Other extensions include Deep Learning strategies such as the Deep Q-Network (DQN) approach presented in \cite{mnih2015human} and the Actor-Critic paradigm \cite{Sutton1998},\cite{mnih2016asynchronous},\cite{lillicrap2015continuous} among many others. 

In this paper, we present an alternative model-based approach, we name Control-Tutored RL (CTRL), where a feedback control strategy designed with only limited or qualitative knowledge of the system dynamics is used to assist the RL algorithm {\em when needed}. For the sake of clarity, we focus on Q-learning (QL) as a RL algorithm and discuss the  resulting {\em control-tutored Q-learning} (CTQL) algorithm showing that it is well apt to deal with continuous or large state spaces while retaining many of the features of a tabular method. We wish to emphasize that our algorithm is complementary to other existing model-based approaches such as \cite{Fathinezhad2016,brunner2017repetitive}. Indeed, in our setting, the control-tutor supports the process of exploring the optimization landscape by suggesting possible actions based on its partial knowledge of the system dynamics. The learning agent can then deploy the policy suggested by the control-tutor whenever it is unable to find a better action by querying the Q-table. With this respect, the control-tutor supports the process of filling in the elements of the Q-table when it is needed, speeding up as a result the convergence of the learning process when compared to that of the Q-learning when used without the tutor. A related but different idea was independently presented in \cite{rathi2020driving}  where RL is mirrored with a Model Predictive Controller (MPC) and a different strategy is used to orchestrate transitions between RL and MPC. To validate our approach, we apply CTQL to solve the classical benchmark problem of stabilizing an inverted pendulum comparing its performance with Q-learning and the feedback strategy used as a tutor. We find that CTQL obtains better performance and convergence than Q-learning or feedback control on their own, solving the stabilization problem even when they are unable to do so by themselves.

\section{Preliminaries} \label{sec:Preliminaries}
Reinforcement learning is an area of machine learning aimed at studying how sub-optimal policies can be computed from data \cite{bertsekas1996neuro,doi:10.1146/annurev-control-053018-023825} in order to solve dynamic programming problems involving  uncertain dynamics \cite{Sutton1998,bertsekas1996neuro}.
By closely following \cite{matni2019self}, we formulate the RL control problem as a constrained optimal control problem of the form:
\begin{align}
    \max_{\pi}& \ \ \mathbb{E}[J_N^{\pi}], \label{stat1.1}\\ 
    \text{s.t.}  & \ \ \begin{aligned}
        x_{k+1} = f_k(x_k,u_k,w_k), \label{stat1.2}
    \end{aligned}\\
    &\ \ \begin{aligned}
    u_k = \pi_k(x_{0:k},u_{0:k-1}),  \label{stat1.3}
    \end{aligned}\\
    &\ \ \begin{aligned}
    x_k \in \mathcal{X} \label{stat1.4}, u_k \in \mathcal{U},     
    \end{aligned}
\end{align}
where $x_k \in \mathcal{X}$ is the state of the system, $u_k \in \mathcal{U}$ is the control input, $w_k \in \mathcal{W}$ is some process noise, $f_k: \mathcal{X} \times \mathcal{U} \times \mathcal{W} \mapsto \mathcal{X}$ is the system vector field, $x_{0:k} = \{x_0, x_1, ..., x_k\}$ is the system state trajectory from step $0$ to step $k$, $u_{0:k} = \{u_0, u_1, ..., u_{k}\}$ is the sequence of control inputs fed to the system from step $0$ to step $k$. 
The objective function $J_N^{\pi}$ in \eqref{stat1.1} is defined as follows:
\begin{equation} \label{eq:objective}
    J_N^{\pi} = r_N(x_{N}) + \sum_{k=1}^{N-1} r_k(x_{k-1},u_{k-1},x_k),
\end{equation}
where $N$ is the time horizon, $r_k:\mathcal{X}^{2} \times \mathcal{U} \mapsto \mathbb{R}$ is some reward function, $\hat r_N:\mathcal{X}\mapsto \mathbb{R}$ is the reward associated to the terminal state $x_N\in\mathcal{X}$ and $\pi = \{\pi_0,\pi_1,...,\pi_{N-1}\}$ is the control policy with $\pi_k:\mathcal{X}^{k+1} \times \mathcal{U}^{k} \mapsto \mathcal{U}$ being a randomized mapping.
The expected value $\mathbb{E}[J_N^{\pi}]$ in \eqref{stat1.1} is taken with respect to the control policy $\pi$ and the random variables $(x_0,w_{0:N})$ assumed to be defined over a common probability space with known and independent distributions.
In many applications of RL to control problems, the system dynamics encoded by $f_k$ is often unknown so the problem \eqref{stat1.1}-\eqref{stat1.4} cannot be directly solved. In this situation, strategies can be used such as the Q-learning algorithm \cite{Watkins1992}. This algorithm seeks to learn an approximation of the reward-to-go function. Such approximation is often carried by a tabular representation, known as the Q-table, which is filled-in online by successive trial and error experiments. In order to succeed, the method requires enough trials so that the learning agent can explore a large enough region of the optimization landscape and hence identify an approximation of the optimal policy.

In what follows, we discuss how a feedback control law with limited knowledge of the plant model can be embedded in the learning process assisting the learning agent in identifying a sub-optimal solution to the control problem in a lesser number of trials.
%%%%%%%%%%%%%%%%%%%%%%%%%%%%%%%%%%%%%%%%%%%%%%%%%%%%%%%%%%%%%%%%%%%%%%%%%%%%%%%%%%%%%%%%%%%%%%%%%%%
\section{The Control-Tutored Reinforcement Learning} \label{sec:RL}
\begin{figure}[t]
\centering
\resizebox{8.4cm}{6cm}{
\begin{tikzpicture}[auto, node distance=1.8cm]
    
    \node [input, name=input] {};
    \node [block, right of=input] (system) {System};
    % \node [block, below of=system, right of=system, xshift=1cm] (RF) {Reward Function};
    \node [block, below of=system, xshift=-3cm , yshift=-2cm, text width=2.1cm] (RL) {Reinforcement Learning};
    \node [block, below of=RL] (CT) {Control Tutor};
    % \node [block, below of=system, left of=system, xshift=-2cm] (policyQ) {\large $\pi^R$};
    % \node [block, below of=system, left of=system, xshift=-4.4cm] (policyT) {\large $\pi^C$};
    % \node [block, below of=policyT, yshift=-0.4cm] (CL) {Control Law};

    %\draw [->] (system.est) -| node [pos=0.25] [name=S] {$x_{k+1},x_{k}$} (RL);
    \draw [-] (system) --  ++(1.4cm,0) coordinate(S1) {};
    \draw [->] (S1) -- node [pos=0.3] [name=A1] [above] {$x_{k}$} ++(0.5cm,0) coordinate(S2) {};
    \draw [-] (S1) |-  ++(0,-4.05cm) coordinate(S2) {};
    \draw [->] (S2) |-  ++(-3.21cm,0) coordinate(S3);
    \draw [->] (S2) |-  (CT);
    % \draw [->] (RF) |-  node [pos=0.1] {$r_{k+1}$} (LU);
    % \draw [->] (S1) |- (Q);
    % \draw [->] (Q) -| (policyQ);
    \draw [<-] (system.west) -- ++(-0.6cm,0) coordinate(p1){};
    \draw [-o] (p1) -- node [pos=0.3] [name=A] [above] {$u_{k}$} ++(-1.1cm,0) coordinate(p2){};
    \draw [-] (p1) |- ++(0,-3.6cm) coordinate(p3);
    \draw [->] (p3) |-  ++(-0.32cm,0) coordinate(p4);
    % \draw [-] (p1) -- ++(0,-0.5cm) coordinate(p2){};
    % \draw [->] (A) |- (RF);
    % \draw [->] (S1) |- (LU);
    % \draw [->] (A) |- (LU);
    % \draw [->] (LU) -- node [pos=0.5] [name=temp] [right] {} (Q);
    % \draw [->] (temp) -| (policyQ);
    % \draw [->] (S1) |- ([xshift=-0.3cm,yshift=-0.27cm]LU.south west) -| (CL);
    % \draw [->] (CL) -- node [pos=0.5] [name=v] [right] {{$v_{k}$}} (policyT);
    \draw [-] (RL) -- node [pos=0.5] [name=A1] [above] {$\pi^R_{k}$} ++ (-2.5cm,0) coordinate(pq1){};
    \draw [-] (pq1) -- ++(0,3.5cm) coordinate(pq2){};
    \draw [-o] (pq2) -- ++(2.2cm,0) coordinate(pq3){};
    \draw [->] (RL.north) -- node [pos=0.5] [name=A2] {$\zeta$} ++ (0,2.6cm) coordinate(pq4){};
    
    \draw [-] (CT) -- node [pos=0.3] [name=A3] [above] {$\pi^C_{k}$} ++(-3.5cm,0) coordinate(pt1){};
    \draw [-] (pt1) -- ++(0,6cm) coordinate(pt2){};
    \draw [-o] (pt2) -- ++(3.2cm,0) coordinate(pt3){};

    \draw [-] (p2) -- (pq3);
    \draw[<-] ($(pt3)+(0.5cm,0.1cm)$) to [bend right]($(pq3)+(0.5cm,-0.1cm)$);
    % \draw[->]  (pq1)  to [out=-150,in=-30] (pt1);
     
\end{tikzpicture}
}
\caption{Schematic of the Control-Tutored Reinforcement Learning (CTRL) algorithm. At each  $k$, the  agent selects its next control action $u_k$ from a given system state $x_k$. This is done by choosing either the control action suggested by the control-tutor policy $\pi_k^C$ or the one suggested by the RL policy $\pi_k^R$. The choice is made in accordance to a boolean variable, $\zeta$, that will be formally defined in Section \ref{sec:CTQL}.} 
\label{fig:CTQL}
\end{figure}
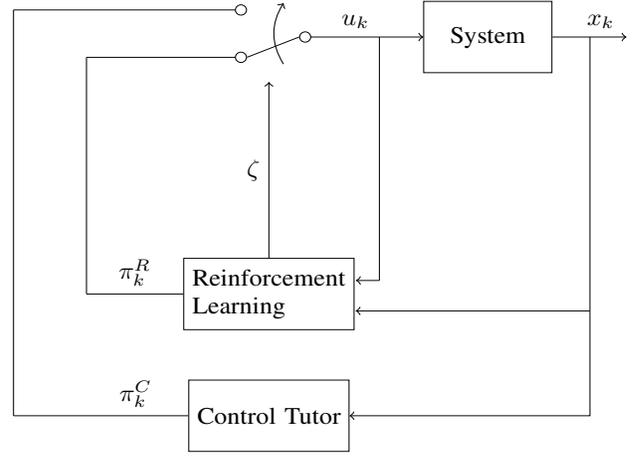

Our approach starts from the observation that in many control problems it is reasonable to assume the knowledge of a mathematical model that partially describes the dynamics of the system. That is, we make the following:
\begin{assumption} \label{assum}
    Only an estimate, say $\hat f_k(x_k,u_k)$, of the system vector field $f_k(x_k,u_k,w_k)$ is available and
    \begin{equation*}
      f_k(x_k,u_k,w_k) = \hat f_k(x_k,u_k) + \delta_k (x_k,u_k,w_k),  
    \end{equation*}
    where $\delta_k (x_k,u_k,w_k)$ is a vector field encompassing all the unknown terms in the dynamics (e.g. higher order terms, unmodelled dynamics etc.)
\end{assumption}
Using $\hat f$ we construct a feedback control strategy (or control-tutor policy)  that the learning agent can use to decide what the next action to take should be. Namely, as shown in Figure \ref{fig:CTQL}, the algorithm will pick either the action suggested by the control-tutor policy, $\pi^C$ or that proposed by the RL algorithm, $\pi^R$, according to which yields the lower expected value of the objective function, $J_N^\pi$. 
\begin{remark}
A key difference between the approach presented here and the one of e.g. \cite{7039601,9140024} is that we do not seek to learn an approximation for $\delta_k (x_k,u_k,w_k)$. Interestingly, as we shall see, the presence of the tutor makes it possible to learn an optimal policy without learning the uncertainty.
\end{remark}
Given the set-up described above, the problem of designing a model-based control-tutor for the RL algorithm (simply termed as the {\em CTRL problem} in what follows) can be stated by adapting the formulation \eqref{stat1.1}-\eqref{stat1.4} as follows:
\begin{align} 
    \max_{\pi} & \ \ \mathbb{E}[J_N^{\pi}], \label{stat2.1}\\ 
    \text{s.t.} & \ \ x_{k+1} = f_k(x_k,u_k,w_k), \label{stat2.2}\\
    & \ \ u_k = \pi_k(x_{0:k},u_{0:k-1},\zeta) \\
    & \ \ \ \ \   = \begin{cases} \label{stat2.3}
\pi^R_k(x_{0:k},u_{0:k-1}) \ \ \ \ \text{if } \zeta \ \text{is true,}\\
\pi^C_k(x_{0:k},u_{0:k-1}) \ \ \ \ \text{otherwise,}
\end{cases}\\
    & \ \ x_k \in \mathcal{X}, u_k \in \mathcal{U}, \label{stat2.4}
\end{align}
where the control policy $\pi_k$ is now selected between $\pi^R_k:\mathcal{X}^{k+1}\times \mathcal{U}^{k} \mapsto \mathcal{U}$, i.e. the randomized mapping defined according to the RL algorithm, and $\pi^C_k:\mathcal{X}^{k+1} \times \mathcal{U}^{k} \mapsto \mathcal{U}$ which is a randomized mapping defined by the feedback control law, and $\zeta$ is a Boolean condition based on the expected value of the reward-to-go.

%%%%%%%%%%%%%%%%%%%%%%%%%%%%%%%%%%%%%%%%%%%%%%%%%%%%%%%%%%%%%%%%%%%%%%%%%%%%%%%%%%%%%%%%%%%%%%%%%%%
\section{Control-Tutored Q-learning (CTQL) implementation} \label{sec:CTQL}
As a representative implementation of CTRL, we develop a control-tutored Q-learning (CTQL) approach by extending the Q-learning approach, briefly described in Sec. \ref{sec:Preliminaries}, to solve a stabilization problem. We make the standard assumption that the policy is time-invariant and that the Markov property holds. 

We select the boolean switching criterion in (\ref{stat2.3})  as
\begin{equation} \label{eq:zeta}
    \zeta = \begin{cases} 
    1 & \text{if}\ \max\limits_{u \in \mathcal{U}}\{Q(x_k,u)\}>0,\\
    0 & \text{otherwise},
    \end{cases}
\end{equation}
%so that the policy becomes:
%\begin{equation} \label{eq:4-6}
%    \pi_k(x_k,\zeta) = 
%    \begin{cases} 
%        \pi^R(x_k)   &    \text{if}\ \zeta \text{\ is true},\\
%        \pi^C(x_k),  &   \text{otherwise}.
%    \end{cases}
%\end{equation}
where $Q(x_k,u)$ is the value stored in the Q-table for the state-action pair $(x_k,u)$ that approximates the reward-to-go from that state. 
In this way at step $k$, given the state $x_k$, the learning agent checks the sign of the entries of the Q-table for all actions $u \in \mathcal{U}$. If at least one of these entries is positive, then the control action $u_k$ is selected according to the classic $\varepsilon$-greedy Q-learning algorithm \cite{Sutton1998}:
\begin{equation} \label{eq:4-5}
  \pi^R(x_k) = 
    \begin{cases} 
        \arg \max\limits_{u\in \mathcal{U}} Q(x_k,u) &   \text{with probability ($1 -\varepsilon$)},\\
        \mathrm{rand}(u)   &   \text{with probability $\varepsilon$}.
    \end{cases}
\end{equation}
Intuitively, this means that the approximate reward-to-go from the current state stored in the Q-table is positive and hence a reward increase is possible by taking this choice.

Otherwise, if the current approximation contained in the Q-table has no positive values, the action is chosen according to the one suggested by the control-tutor via the policy  $\pi^C(x_k)$ which we also choose as a $\varepsilon$-greedy policy of the form:
\begin{equation} \label{eq:4-4}
  \pi^C(x_k) = 
    \begin{cases} 
        \arg \min\limits_{ u \in \mathcal{U}}\|v(x_k) - u\| &   \text{with probability ($1 -\varepsilon$)},\\
        \mathrm{rand}(u)   &   \text{with probability $\varepsilon$},
    \end{cases}
\end{equation}
where $v$ is the control input generated by a feedback controller designed using the model estimate $\hat f$ available to the control-tutor. As such input does not necessarily belong to $\mathcal{U}$, the policy function $\pi^C(x_k)$ selects the action $u\in \mathcal{U}$ which is closest to $v(x_k)$.

Once the action is selected from either $\pi^R(x_k)$ or $\pi^C(x_k)$, the corresponding expected reward is then computed and used to update the Q-table. The pseudo-code of the CTQL algorithm is given in Algorithm \ref{alg:CTQL}. 

Note that to preserve the spirit of the Q-learning algorithm  both policies $\pi^C(x_k)$ and $\pi^R(x_k)$ contain some degree of randomness to favour exploration. We are currently studying whether this guarantees that, when implemented, the policy selection function of the CTQL is still within the scope of the probabilistic proof of convergence available for the Q-learning algorithm and described in \cite{Watkins1992,bertsekas1996neuro}.

\begin{algorithm}[htbp]
  \caption{control-tutored Q-learning}
  \label{alg:CTQL}
  \begin{algorithmic}
    \STATE {Initialize $Q(x,u) = 0,\forall x \in \mathcal{X},u \in \mathcal{U}$}
      \STATE{Detect intial state $x_0$}
      \FOR{$k = 0$ to $N$}
      \STATE{Compute $\zeta$ via \eqref{eq:zeta}}
        \IF{$\zeta$}
           \STATE{$u_k \gets\pi^R(x_k)$}
       \ELSE
           \STATE{$u_k \gets\pi^C(x_k)$}
       \ENDIF{\bf{end if}}
       \STATE {Observe and store $x_{k+1}$ and $r_{k}$}
       \STATE {$Q(x_{k},u_{k})\gets(1-\alpha)Q(x_{k},u_{k}) + \alpha[r_{k} + $}
       \STATE {$\ \ \ \ \ \ \ \ \ \ \ \ \ \ \ \ \ \ \ \ + \gamma \max\limits_{u \in U}Q(x_{k+1},u)]$}
      \ENDFOR{\bf{end for}}
  \end{algorithmic}
\end{algorithm}

To illustrate the viability and effectiveness of CTQL, we apply it to solve the problem of stabilizing an inverted pendulum and discuss its performance by comparing it to a traditional (untutored) Q-learning approach.
%%%%%%%%%%%%%%%%%%%%%%%%%%%%%%%%%%%%%%%%%%%%%%%%%%%%%%%%%%%%%%%%%%%%%%%%%%%%%%%%%%%%%%%%%%%%%%%%%%%
\section{Application to the Inverted Pendulum}
We consider the problem of stabilizing the physical pendulum provided by the OpenAI Gym framework \cite{brockman2016openai,GYM} in its inverted position. 
\subsection{Problem Formulation}
To achieve the control objective we define the reward function in (\ref{eq:objective}) as:
\begin{equation} \label{eq:reward}
    r_{k}(x_{k},x_{k-1}) = - \left[V(x_{k}) - V(x_{k-1})\right] + \rho(x_k),
\end{equation}
where $x_k=[x_{k,1},\ x_{k,2}]$ is the system state with $x_{k,1}$ and $x_{k,2}$ being the angular position and angular velocity of the pendulum respectively, meanwhile $V:\mathcal{X} \mapsto \mathbb{R}$ is a scalar quadratic function defined as:
\begin{equation}
    V(x_k) = k_1x_{1,k}^2 + k_2x_{2,k}^2,
\end{equation}
$\rho(x_k)$ is an additional term of the reward function that accounts for a positive prize $p$ that the agent receives if the pendulum falls in a sufficiently small neighborhood of the upward position. Specifically, such a term is defined as:
\begin{equation} \label{eq:rhat}
    \rho(x_k) = 
        \begin{cases}
            p \ \ \ \text{ if } x_{1,          k} \in [-\epsilon, \epsilon],\\
            0 \ \ \ \text{ otherwise, }
        \end{cases}
\end{equation}
Also, we let $r_N(x_N) = \rho(x_N)$ in \eqref{eq:objective}.
\begin{remark}
Substituting \eqref{eq:reward} in \eqref{eq:objective} we get, in the expression for $J_N^{\pi}$, the sum -$\sum_{k=0}^{N-1}V(x_{k+1}) - V(x_k) = -(V(x_N)-V(x_0))$. In the special case where the function $V(\cdot)$ is the Lyapunov function for the system, maximizing this term implies finding a solution to the optimization problem in \eqref{stat1.1}-\eqref{stat1.4} that minimizes the $N$-step derivative of the Lyapunov function. We note how, in the example described in this section, the CTQL is able to stabilize the upward equilibrium of the inverted pendulum even if $V(\cdot)$ is not chosen a Lyapunov function for the system\footnote{We leave for future research the problem of finding analytical conditions on $V(\cdot)$ that make it a viable choice for the CTQL algorithm to solve a given problem.}.  
\end{remark}
We start by implementing the classical Q-learning algorithm. The state space is defined as $\mathcal X:= \mathcal{D}\times\mathcal{G}$, where $\mathcal{D}$ is the set of angular positions and $\mathcal{G}$ is the set of angular speeds. The action space $\mathcal{U}$ is the set of possible values of the control input $u_k$. In our implementation, the sets $\mathcal{D},\mathcal{G}$ and $\mathcal{U}$ are obtained by discretizing the continuous state space and control input of the pendulum as described in the Appendix. 
The results of the Q-learning implementation on its own are shown in Fig. \ref{fig:Rewards} and Fig. \ref{fig:TQL}. 
\input{Figures/Rewards}
\input{Figures/QL_Traj}
%%%%%%%%%%%%%%%%%%%%%%%%%%%%%%%%%%%%%%%%%%%%%%%%%%%%%%%%%%%%%%%%%%%%%%%%%%%%%%%%%%%%%%%%%%%%%%%%%%%
The design of the tutor control law requires some model of the expected dynamics. 
We assume that only the linearized version of the inverted pendulum model around the upward position $x_k = [0,0]^T$ is available to the control-tutor. Specifically we assume the following model of the pendulum dynamics is available to the tutor:
\begin{equation}\label{eq:4-3}
x_{k+1} = \hat f(x_k,v_k)= Ax_k + Bv_k, 
\end{equation} 
with $B = [0,\ dt\frac{1}{I}]^T$ and $A$ defined as:
\begin{equation}
    A = \begin{bmatrix}
    0 & 1+dt \\
    g\frac{l}{2I}dt & 1
    \end{bmatrix},
\end{equation}
where $dt$ is the sampling time, $I = m\frac{l^2}{3}$ is the inertia of the homogeneous rod, $g$ is the gravitational constant and $l$ and $m$ are the length and the mass of the rod respectively. All parameters values chosen for this study can be found in the Appendix.

The control-tutor is then designed as the state-feedback control input:
\begin{equation} \label{control input}
v_k=-Kx_k,
\end{equation}
with the control gains selected as $K = [5.83, 1.83]^T$ in order to render the origin of the linearized system a stable node with a settling time of $10s$ (corresponding to about 200 steps in our discretization).
The control input $v_k$ defined in \eqref{control input} is then used in \eqref{eq:4-4} to obtain the control-tutor policy. The results of the CTQL implementation are shown in Fig. \ref{fig:Rewards} and Fig. \ref{fig:TCTQL}.
\input{Figures/CTQL_Traj}
%%%%%%%%%%%%%%%%%%%%%%%%%%%%%%%%%%%%%%%%%%%%%%%%%%%%%%%%%%%%%%%%%%%%%%%%%%%%%%%%%%%%%%%%%%%%%%%%%%%
\subsection{Comparison between CTQL and Q-learning}
The numerical validation is carried out on the Pendulum-V0 environment \cite{GYM}. We define the training sessions set $\mathcal{S} = \{1,...,S\}$, the episodes set $\mathcal{E} = \{1,...,E\}$ and the simulation set $\mathcal{N} = \{1,...,N\}$. Each episode corresponds to a simulation of $N$ steps of the pendulum starting from the initial condition on the stable downward position. More details on the simulation parameters can be found in the Appendix.

The numerical results are used to evaluate both the data efficiency of the learning process and the control performance. 
In what follows we use the superscript $e$ to denote the episode at which the variables are evaluated. We denote as ${\mathcal{M}}$ the set of consecutive episodes where the learning agent is able to maintain the pendulum position $x_{1,k}$, and velocity, $x_{2,k}$ in a ball of radius $\epsilon$ in the time window $k \in [N/2, N]$. We then assume the learning phase ends if the following condition is satisfied: 
\begin{equation} \label{eq:terminal}
    \sum_{e \in {\mathcal{M}}}\big(r_N^e(x_N) +  \sum_{k=N/2}^{N-1} {r_k^{e}(x_k,x_{k-1})}\big) \geq M\frac{N}{2}p,
\end{equation}
where $M$ is the minimum number of successful consecutive episodes we set as a satisfactory threshold (see the Appendix for the value we chose in our implementation).

\subsubsection{Learning performance}

to evaluate the sample efficiency and to make a proper comparison between the strategies proposed, we define the following metrics:
\begin{equation} \label{l1}
    \lambda_1 = \max_{e\in\mathcal{E}}J_N^{\pi,e},
\end{equation}
\begin{equation} \label{l2}
    \lambda_2 = \frac{1}{E}\sum_{e\in\mathcal{E}} J_N^{\pi,e},
\end{equation}
where $\lambda_1$ is the maximum value of the objective function $J_N^{\pi}$, $\lambda_2$ is the average value of the objective $J_N^{\pi}$ with respect to the total number of episodes. 

We also define $\lambda_3$ as the number of episodes after which the terminal condition \eqref{eq:terminal} is satisfied.

The metrics defined above are averaged over $S=5$ runs of the algorithms with $E=6000$ episodes each for both Q-learning on its own and CTQL. The results are summarized in Table \ref{tab:learning} showing that both achieve a comparable value of $\lambda_1$ but CTQL guarantees a better average value of the objective function across all episodes (metric $\lambda_2$), fulfilling the terminal condition \eqref{eq:terminal} after a notably smaller number of episodes (metric $\lambda_3$).
\begin{table}[htbp]
\caption{Data efficiency and learning performance comparison between Q-learning and CTQL}
\label{tab:learning}
\begin{center}
\begin{tabular}{||c||c||c||c||}
\hline
 & $\lambda_1$ & $\lambda_2$ & $\lambda_3$\\
\hline
QL & 1552 & 763 & 2730\\
\hline
CTQL & 1521 & 1240 & 417 \\
\hline
\end{tabular}
\end{center}
\end{table}
To further compare the two strategies we show in Fig. \ref{fig:Tut} the number of times per episode in which the action suggested by the control-tutor is taken by the learning agent. We observe that the control-tutor policy is most deployed by the agent during the initial episodes with the number of implemented actions coming from the tutor decreasing as the agent converges towards a suitable control strategy fulfilling the control goal. 
\input{Figures/Tutor}
%%%%%%%%%%%%%%%%%%%%%%%%%%%%%%%%%%%%%%%%%%%%%%%%%%%%%%%%%%%%%%%%%%%%%%%%%%%%%%%%%%%%%%%%%%%%%%%%%%%
\subsubsection{Control performance}
finally, we define a set of metrics to evaluate and compare the performance of the controller obtained at the end of the learning stage using the CTQL and QL algorithms. Specifically, the control metrics are defined as follows:

\begin{itemize}
    \item the settling time $\eta_1$ defined as the first step such that
    \begin{equation} \label{eq:eta1}
     \|x_{k}\|\leq 0.05  \ \ \forall k \geq \eta_1;
\end{equation}
\item the average value $\eta_2$ of the mismatch between the state and the target equilibrium (i.e. the origin) over the last $\Gamma=100$ steps:
\begin{equation} \label{eq:eta2}
    \eta_2= \frac{1}{\Gamma}\sum^N_{k=N-\Gamma}\|x_{k}\|;
\end{equation}
\item the value of the objective computed over the trajectory as defined in (\ref{eq:objective}), i.e. $\eta_3:= J_N^{\pi}$.
\end{itemize}

The control performances metrics were evaluated by running the controller after the end of the learning phase to swing up the pendulum from its stable downward position. The results are summarized in Table \ref{tab:control} where we see that the control performance of controller trained using the CTQL is comparable to that of the controller obtained by running the QL algoritm; the notable difference being therefore the much quicker learning times of CTQL as summarized in Table \ref{tab:learning}.

\begin{table}[htbp]
\caption{Control performance comparison between the controller trained using Q-learning and that trained using a Control-Tutored Q-learning approach}
\label{tab:control}
\begin{center}
\begin{tabular}{||c||c||c||c||}
\hline
 & $\eta_1$ & $\eta_2$ & $\eta_3$\\
\hline
QL & 125 & 0.00968 &  1384\\
\hline
CTQL & 129 & 0.00901 & 1359 \\
\hline
\end{tabular}
\end{center}
\end{table}

%%%%%%%%%%%%%%%%%%%%%%%%%%%%%%%%%%%%%%%%%%%%%%%%%%%%%%%%%%%%%%%%%%%%%%%%%%%%%%%%%%%%%%%%%%%%%%%%%%%
\section{Conclusions}
We introduced an extension of reinforcement learning where the policy selection function is enhanced by means of a control-tutor that, using a feedback control law with limited knowledge of the system dynamics, is able to support the exploration of the optimization landscape guaranteeing better convergence and shorter learning times. To illustrate the effectiveness of the approach, we tutored the Q-learning algorithm via a state feedback controller to solve the classical benchmark problem of stabilizing an inverted pendulum. In our experiments, the tutor only had access to a linearized model of the pendulum about its inverted position. In this situation, we showed that our CTQL strategy is able to swing up the pendulum and globally stabilize its inverted position with the learning process converging towards a viable control solution after a much shorter number of episodes than the QL when deployed on its own. Ongoing work is focused on refining this approach with the aim of giving theoretical guarantees, obtaining a better understanding of its advantages and limitations for future applications.

\addtolength{\textheight}{-1.5cm}   % This command serves to balance the column lengths
                                  % on the last page of the document manually. It shortens
                                  % the textheight of the last page by a suitable amount.
                                  % This command does not take effect until the next page
                                  % so it should come on the page before the last. Make
                                  % sure that you do not shorten the textheight too much.

%%%%%%%%%%%%%%%%%%%%%%%%%%%%%%%%%%%%%%%%%%%%%%%%%%%%%%%%%%%%%%%%%%%%%%%%%%%%%%%%

%%%%%%%%%%%%%%%%%%%%%%%%%%%%%%%%%%%%%%%%%%%%%%%%%%%%%%%%%%%%%%%%%%%%%%%%%%%%%%%%

%%%%%%%%%%%%%%%%%%%%%%%%%%%%%%%%%%%%%%%%%%%%%%%%%%%%%%%%%%%%%%%%%%%%%%%%%%%%%%%%
\section*{Acknowledgements}
The authors wish to acknowledge the contributions of Prof Pietro De Lellis, University of Naples Federico II, and Ms. Fabrizia Auletta, University of Bristol and Macquarie University, to an earlier version of the CTQL algorithm and its application to the solution of an herding problem that was reported earlier in \cite{de2019control}.

\section*{Appendix}
We give here all the parameters that were used for the numerical simulations reported in the paper.
We define the training set composed of $S=5$ learning sessions each of them composed of $E=6000$ episodes, each corresponding to a simulation of $N=400$ steps of the pendulum starting from its stable downward position.
%in page 157 of \cite{Sutton1998} 
The Q-learning update rule is implemented as showed in Algorithm \ref{alg:CTQL} with parameters set to  $\alpha = 0.8$ and  $\gamma = 0.97$ while the randomness parameter in the policies $\pi^R(\cdot)$ and $\pi^C(\cdot)$ is set to $\varepsilon = 0.02$. 
The parameters of the reward function are set to $k_1 = 1$, $k_2 = 0.1$ while the price for reaching a state in the neighborhood of the upward position is set as $p=5$ and it is given in the region $x_1\in[-\epsilon,\epsilon]$ with $\epsilon = 0.05$.
The terminal condition is met when the agent stabilizes the pendulum for at least half of the $N=400$ steps for $M=20$ consecutive episodes. The state space $\mathcal{X}$ discretization is non-uniform and chosen as follows. The angular position is discretized in 16 equally spaced values when $x_{1} \in [-\pi,-\frac{\pi}{9}] \cup (\frac{\pi}{9},\pi]$ while it is split in 14 discrete values for  $x_{1} \in (-\frac{\pi}{9},-\frac{\pi}{36}] \cup [\frac{\pi}{36},\frac{\pi}{9})$ and 10 equally spaced values for $x_{1} \in (-\frac{\pi}{36},\frac{\pi}{36})$. The angular velocity is discretized in 20 equally spaced values when $x_{2} \in [-8,-1) \cup (1,8]$ and 18 equally spaced values when $x_{2} \in [-1,1]$. The control action space  $\mathcal{U}$ discretization is discretized as follows: 18 equally spaced values for $u \in [-2,-0.2] \cup [0.2,2]$ and 8 equally spaced values for $u \in (-0.2,0.2)$.

\begin{thebibliography}{10}
\providecommand{\url}[1]{#1}
\csname url@rmstyle\endcsname
\providecommand{\newblock}{\relax}
\providecommand{\bibinfo}[2]{#2}
\providecommand\BIBentrySTDinterwordspacing{\spaceskip=0pt\relax}
\providecommand\BIBentryALTinterwordstretchfactor{4}
\providecommand\BIBentryALTinterwordspacing{\spaceskip=\fontdimen2\font plus
\BIBentryALTinterwordstretchfactor\fontdimen3\font minus
  \fontdimen4\font\relax}
\providecommand\BIBforeignlanguage[2]{{%
\expandafter\ifx\csname l@#1\endcsname\relax
\typeout{** WARNING: IEEEtran.bst: No hyphenation pattern has been}%
\typeout{** loaded for the language `#1'. Using the pattern for}%
\typeout{** the default language instead.}%
\else
\language=\csname l@#1\endcsname
\fi
#2}}

\bibitem{Sutton1998}
R.~S. Sutton and A.~G. Barto, \emph{Reinforcement Learning: An Introduction},
  2nd~ed.\hskip 1em plus 0.5em minus 0.4em\relax The MIT Press, 2018.

\bibitem{bertsekas1996neuro}
D.~P. Bertsekas and J.~N. Tsitsiklis, \emph{Neuro-dynamic programming}.\hskip
  1em plus 0.5em minus 0.4em\relax Athena Scientific, 1996.

\bibitem{kober2013reinforcement}
J.~Kober, J.~A. Bagnell, and J.~Peters, ``Reinforcement learning in robotics: A
  survey,'' \emph{The International Journal of Robotics Research}, vol.~32,
  no.~11, pp. 1238--1274, 2013.

\bibitem{garcia2015comprehensive}
J.~Garc{\i}a and F.~Fern{\'a}ndez, ``A comprehensive survey on safe
  reinforcement learning,'' \emph{Journal of Machine Learning Research},
  vol.~16, no.~1, pp. 1437--1480, 2015.

\bibitem{cheng2019end}
R.~Cheng, G.~Orosz, R.~M. Murray, and J.~W. Burdick, ``End-to-end safe
  reinforcement learning through barrier functions for safety-critical
  continuous control tasks,'' \emph{preprint available from arXiv:1903.08792},
  2019.

\bibitem{berkenkamp2017safe}
F.~Berkenkamp, M.~Turchetta, A.~Schoellig, and A.~Krause, ``Safe model-based
  reinforcement learning with stability guarantees,'' \emph{Advances in neural
  information processing systems}, vol.~30, pp. 908--918, 2017.

\bibitem{gu2016continuous}
S.~Gu, T.~Lillicrap, I.~Sutskever, and S.~Levine, ``Continuous deep q-learning
  with model-based acceleration,'' \emph{Proc. of the International Conference
  on Machine Learning}, pp. 2829--2838, 2016.

\bibitem{deisenroth2011pilco}
M.~Deisenroth and C.~E. Rasmussen, ``Pilco: A model-based and data-efficient
  approach to policy search,'' \emph{Proc. of the International Conference on
  Machine Learning}, pp. 465--472, 2011.

\bibitem{rosolia2017learning}
U.~Rosolia and F.~Borrelli, ``Learning model predictive control for iterative
  tasks. a data-driven control framework,'' \emph{IEEE Transactions on
  Automatic Control}, vol.~63, no.~7, pp. 1883--1896, 2017.

\bibitem{rathi2020driving}
M.~Rathi, P.~Ferraro, and G.~Russo, ``Driving reinforcement learning with
  models,'' \emph{Proceedings of SAI Intelligent Systems Conference}, pp.
  70--85, 2020.

\bibitem{mnih2015human}
V.~Mnih, K.~Kavukcuoglu, D.~Silver, A.~A. Rusu, J.~Veness, M.~G. Bellemare,
  A.~Graves, M.~Riedmiller, A.~K. Fidjeland, G.~Ostrovski, \emph{et~al.},
  ``Human-level control through deep reinforcement learning,'' \emph{Nature},
  vol. 518, no. 7540, p. 529, 2015.

\bibitem{mnih2016asynchronous}
V.~Mnih, A.~P. Badia, M.~Mirza, A.~Graves, T.~Lillicrap, T.~Harley, D.~Silver,
  and K.~Kavukcuoglu, ``Asynchronous methods for deep reinforcement learning,''
  in \emph{Proc. of International conference on machine learning}, 2016, pp.
  1928--1937.

\bibitem{lillicrap2015continuous}
T.~P. Lillicrap, J.~J. Hunt, A.~Pritzel, N.~Heess, T.~Erez, Y.~Tassa,
  D.~Silver, and D.~Wierstra, ``Continuous control with deep reinforcement
  learning,'' \emph{preprint available from arXiv:1509.02971v6}, 2019.

\bibitem{Fathinezhad2016}
F.~Fathinezhad, V.~Derhami, and M.~Rezaeian, ``Supervised fuzzy reinforcement
  learning for robot navigation,'' \emph{Applied Soft Computing}, vol.~40, pp.
  33 -- 41, 2016.

\bibitem{brunner2017repetitive}
M.~Brunner, U.~Rosolia, J.~Gonzales, and F.~Borrelli, ``Repetitive learning
  model predictive control: An autonomous racing example,'' \emph{Proc. of the
  IEEE Conference on Decision and Control}, pp. 2545--2550, 2017.

\bibitem{doi:10.1146/annurev-control-053018-023825}
B.~Recht, ``A tour of reinforcement learning: The view from continuous
  control,'' \emph{Annual Review of Control, Robotics, and Autonomous Systems},
  vol.~2, no.~1, pp. 253--279, 2019.

\bibitem{matni2019self}
N.~Matni, A.~Proutiere, A.~Rantzer, and S.~Tu, ``From self-tuning regulators to
  reinforcement learning and back again,'' \emph{IEEE 58th Conference on
  Decision and Control (CDC)}, pp. 3724--3740, 2019.

\bibitem{Watkins1992}
C.~J. C.~H. Watkins and P.~Dayan, ``Q-learning,'' \emph{Machine Learning},
  vol.~8, pp. 279--292, 1992.

\bibitem{7039601}
A.~K. {Akametalu}, J.~F. {Fisac}, J.~H. {Gillula}, S.~{Kaynama}, M.~N.
  {Zeilinger}, and C.~J. {Tomlin}, ``Reachability-based safe learning with
  gaussian processes,'' in \emph{53rd IEEE Conference on Decision and Control},
  2014, pp. 1424--1431.

\bibitem{9140024}
M.~{Greeff} and A.~P. {Schoellig}, ``Exploiting differential flatness for
  robust learning-based tracking control using gaussian processes,'' \emph{IEEE
  Control Systems Letters}, vol.~5, no.~4, pp. 1121--1126, 2021.

\bibitem{brockman2016openai}
G.~Brockman, V.~Cheung, L.~Pettersson, J.~Schneider, J.~Schulman, J.~Tang, and
  W.~Zaremba, ``Open{AI} {Gym},'' \emph{arXiv preprint arXiv:1606.01540}, 2016.

\bibitem{GYM}
\BIBentryALTinterwordspacing
{Open{AI}}, \emph{{Open{AI} {Gym} Pendulum-v0}}, 2019. [Online]. Available:
  \url{https://github.com/openai/gym/blob/master/gym/envs/classic_control/pendulum.py}
\BIBentrySTDinterwordspacing

\bibitem{de2019control}
F.~De~Lellis, F.~Auletta, G.~Russo, P.~De~Lellis, and M.~di~Bernardo,
  ``Control-tutored reinforcement learning,'' \emph{arXiv preprint
  arXiv:1912.06085}, 2019.

\end{thebibliography}
\end{document}